\documentclass[11pt]{amsart}
\usepackage{amssymb, amscd, graphicx}

\newcommand{\hol}{\mathrm{hol}}

\newcommand{\bC}{\mathbb{C}}

\newcommand{\bP}{\mathbb{P}}

\newcommand{\bR}{\mathbb{R}}
\newcommand{\bZ}{\mathbb{Z}}

\newcommand{\cA}{\mathcal{A}}

\newcommand{\cD}{\mathcal{D}}

\newcommand{\cL}{\mathcal{L}}

\newcommand{\cO}{\mathcal{O}}

\newcommand{\cY}{\mathcal{Y}}

\newcommand{\cT}{\mathcal {T}}

\newcommand{\dist}{\mathrm{dist}}

\theoremstyle{definition}
\newtheorem{lemm}{Lemma}
\newtheorem{prop}{Proposition}
\newtheorem{defi}{Definition}

\newtheorem*{theo*}{Theorem}
\newtheorem*{prop*}{Proposition}
\newtheorem*{rema}{Remark}
\newtheorem*{ackn}{Acknowledgements}

\begin{document}

\title{Homological Mirror Symmetry is T-duality for $\mathbb P^n$}
\author{Bohan Fang}
\address{Department of Mathematics, Northwestern University,
2033 Sheridan Road, Evanston, IL  60208}
\email{b-fang@math.northwestern.edu}

\begin{abstract}
In this paper, we apply the idea of T-duality to
projective spaces. From a connection on a line bundle on $\mathbb
P^n$, a Lagrangian in the mirror Landau-Ginzburg model is
constructed. Under this correspondence, the full strong
exceptional collection $\mathcal O_{\mathbb
P^n}(-n-1),\dots,\mathcal O_{\mathbb P^n}(-1)$ is mapped to
standard Lagrangians in the sense of \cite{nz}. Passing to
constructible sheaves, we explicitly compute the quiver structure
of these Lagrangians, and find that they match the quiver
structure of this exceptional collection of $\mathbb P^n$. In this
way, T-duality provides quasi-equivalence of the Fukaya category
generated by these Lagrangians and the category of coherent sheaves
on $\mathbb P^n$, which is a kind of homological mirror symmetry.
\end{abstract}
\maketitle

\section{Introduction}

Mirror symmetry was first observed for Calabi-Yau manifolds. It
has been extended to Fano cases by considering Landau-Ginzburg
models as mirrors to Fano varieties \cite{hv}. A Landau-Ginzburg
model is a noncompact manifold equipped with a complex-valued
holomorphic function called the superpotential. In this paper, we
study the approach of homological mirror symmetry proposed by
Kontsevich \cite{k2}. Kontsevich suggests to investigate
homological mirror symmetry as the equivalence of the Fukaya
category on the A-model side and the category of coherent sheaves
on its mirror B-model for Calabi-Yau manifolds. Kontsevich
\cite{k} and Hori-Iqbal-Vafa\cite{hiv} discuss homological mirror
symmetry in the case of Fano manifolds. Following this line,
Auroux, Katzarkov and Orlov prove homological mirror symmetry for
weighted projective planes (and their non-commutative
deformations) \cite{ako} and Del Pezzo surfaces \cite{ako2}.
Abouzaid proves the case of all smooth projective toric varieties
using tropical geometry \cite{ab} \cite{ab2}. Bondal and Ruan also
announce a result for weighted projective spaces \cite{b}.

Strominger, Yau and Zaslow \cite{syz} conjecture that mirror
symmetry is a manifestation of T-duality on a special Lagrangian
torus fibration. In the case of toric Fano varieties \cite{lv}
\cite{h}, the moment map produces a fibration by Lagrangian tori.
Auroux \cite{au} discusses the relation between the SYZ conjecture
and Fano varieties. This paper deals with homological mirror
symmetry for the projective space $\mathbb P^n$, using the
philosophy of T-duality. We will apply this T-duality to any
holomorphic line bundle over $\mathbb P^n$, and obtain a
Lagrangian in the Landau-Ginzburg mirror. We \emph{define} $n+1$
Lagrangians $\mathcal L(-1),\dots,\mathcal L(-n-1)$ in the mirror
which arise via T-duality of the collection $\mathcal
O(-1),\dots,\mathcal O(-n-1)$. We roughly state our theorem here.
\begin{theo*}
The derived Fukaya category containing $\mathcal
L(-1),\dots,\mathcal L(-n-1)$ is equivalent to the derived
category of coherent sheaves on $\mathbb P^n$.
\end{theo*}

The definition of this Fukaya category will be specified later. We
remark that the Lagrangians we are considering are submanifolds in
$(\bC^*)^n$ considered as the cotangent bundle of $(S^1)^n$. While $(\bC^*)^n$ together with a certain superpotential $W$ is the Landau-Ginzburg mirror of $\bP^n$ in the sense of \cite{k} and \cite{hiv}, we do not explicitly consider the superpotential $W$ here.
This differs from the case in \cite{ako} and
\cite{ako2}, in which the authors consider the Fukaya-Seidel
category consisting of vanishing cycles in a generic fiber of the
superpotential. However, the Lagrangians in this paper are somehow like the
vanishing thimbles. Under some tentative calculation, the images of the Lagrangiangs $\mathcal
L(k)$ under the superpotential $W$ are not horizontal half-lines
going from critical values to the positive infinity. Instead, they
are ``thickened" rays. Thus the Fukaya category we are using
differs from the original idea of \cite{hiv}. We hope a
renormalization process suggested by \cite{hv} \cite{au} will
remedy this problem in the future.

\begin{rema}
When computing this Fukaya category, we pass to the dg category of constructible sheaves by the results of \cite{nz} and \cite{n}. The constructible sheaves coming from $\cL(-1),\dots, \cL(-n-1)$ are constructible with respect to a particular stratification. This stratification coincides with Bondal \cite{b} and \cite{b2} in the case of projective spaces. The situation in more general cases and the comparison with Bondal's results will be discussed in \cite{fltz}.

The use of an explicit exceptional collection of $\mathbb P^n$ is basically for convenience only. T-duality produces an object for any line bundle, as does the unique functor defined to agree with T-duality on a generating exceptional collection. The maps of objects do not a priori coincide. This issue is addressed in \cite{fltz}, and they turn out to be the same. The method of this paper is extended to treat all projective toric varieties in \cite{fltz}, thus recovering Abouzaid's result \cite{ab2}.
\end{rema}

\begin{ackn}
I would like to thank my advisor, Eric Zaslow, for
showing me the idea of T-duality on projective spaces, and for
valuable discussions and encouragement throughout this project.
I would also like to thank Chiu-Chu Liu and David Treumann, for the discussions when writing a joint paper \cite{fltz}, which lead to many important modifications and simplifications of this one.
\end{ackn}

\section{Mirror of projective spaces}

In this section, we describe the mirror of the projective space
$\mathbb P^n$ over $\mathbb C$. Roughly speaking, the mirror of
$\mathbb P^n$ is a Landau-Ginzburg model, i.e. a non-compact
manifold with a superpotential \cite{k} \cite{hv}. The SYZ
conjecture \cite{syz} suggests to construct the mirror of $\mathbb
P^n$ via T-duality on the torus fibration over the moment
polytope, as argued in \cite{lv}. The superpotential is
conjecturally given by Fukaya-Oh-Ohta-Ono's $m_0$ obstruction
\cite{fooo}. This has been investigated by Cho and Oh \cite{co}.
Following this idea, we give a description of the mirror of
$\mathbb P^n$ as a complexified moduli space of special
Lagrangians. The material in this section is from \cite{au}.

On $\mathbb P^n$ with its Fubini-Study metric and corresponding
symplectic form $\omega$, the torus $T^n$ acts via
$$(\theta_1,\dots,\theta_n)\cdot(z_0:z_1:\dots:z_n)=(z_0:e^{2\pi i \theta_1}z_1:e^{2\pi i
\theta_2}z_2:\dots:e^{2\pi i \theta_n}z_n).$$ This action produces
a moment map $\phi: \mathbb P^n\rightarrow \Delta$ given by
$$ (z_0:z_1:\dots:z_n)\mapsto (\frac{|z_1|^2}{\sum\limits_{i=0}^n
|z_i|^2},\frac{|z_2|^2}{\sum\limits_{i=0}^n
|z_n|^2},\dots,\frac{|z_n|^2}{\sum\limits_{i=0}^n |z_i|^2}),$$
where $$\Delta=\left\{(x_1,x_2,\dots,x_n):x_i\ge 0,\
\sum\limits_{i=1}^{n}x_i\le 1\right\}.$$ This moment map is a
torus fibration outside of the boundary.

Let $D$ be the boundary divisor of $\phi$, i.e. by
$D=\phi^{-1}(\partial \Delta)=\{(z_0:z_1:\dots:z_n)\in\mathbb
P^n|z_0\cdots z_n=0\}$. The holomorphic $n$-form $\Omega$ on
$\mathbb P^n\backslash D$ is $d\log z_1\wedge\dots \wedge d\log
z_n$ in coordinates $(1:z_1:\dots:z_n)$. By a special Lagrangian
we mean a Lagrangian submanifold $L$ with a constant phase
$\vartheta$, i.e. $\mathrm{Im} (e^{-i\vartheta}\Omega)|_L=0$. It
is easy to see that the moment map $\phi$ defines a $T^n$-orbit
fibration on $\mathbb P^n\backslash D$. Because of the following
lemma, this is a special Lagrangian fibration.
\begin{prop}
The $T^n$-orbits in $\mathbb P^n\backslash D$ are special
Lagrangians.
\end{prop}
\proof It is a classical fact that any $T^n$-orbit is a
Lagrangian. Notice the holomorphic form $\Omega$ on $\mathbb
P^n\backslash D$ is $d\log z_1\wedge\dots\wedge d \log z_n$. We
see the restriction of $\Omega$ on the orbit gives rise to phase
$n\pi/2$.\qed

On $\mathbb P^n$, define a polar coordinate system $(r,\theta)$ by
$z_k=r_k e^{i\theta_k}$ at the point $(1:z_1:\dots:z_n)$. The
coordinate $r$ is the coordinate on the base of the fibration
$\phi$, while $\theta$ lives in the fiber.

The mirror $M$ is constructed as the moduli space of the torus
fibers together with flat connections on them.
\begin{defi}
The complexified moduli space $M$ consists of pairs $(L,\nabla)$.
Here $L$ is a $T^n$-orbit in $\mathbb P^n\backslash D$, and
$\nabla$ is a flat $U(1)$ connection on the trivial line bundle
over $L$ up to gauge equivalence.
\end{defi}
We omit the details of the complex and the symplectic structures
on $M$. The result is stated here without any explanation. Any
$T^n$-orbit in $\mathbb P^n\backslash D$ has to be a fiber of
$\phi$. We write a fiber
$$L(r):=S^1(r_1)\times\dots\times
S^1(r_n)=\{(1:z_1:\dots:z_n), |z_i|=r_i\},$$ for
$r=(r_1,\dots,r_n)\in (\mathbb R^+)^n$. It is obvious that
$L(r_1,\dots,r_n)$ is mapped to a point in $\Delta$ by the moment
map $\phi$. The flat connection $\nabla$ on $L$ can be represented
as $\sum \gamma_i d \theta_i$. Therefore we can give a coordinate
system to $M$ by $(r,\gamma)$ where $r=(r_1,\dots,r_n)$ and
$\gamma=(\gamma_1,\dots,\gamma_n)$. Notice $\gamma_i$
takes value in $S^1=\bR/\bZ$.
Since $\mathbb P^n\backslash D$ and $M$ are dual fibrations on the
same base, they share the same coordinate $r$ on the base. With
these coordinates $(r,\gamma)$ in hand, the symplectic structure
on $M$ is defined as following.
\begin{prop}
The symplectic structure on $M$ is given by
$$\omega^\vee=(2\pi)^n\sum d \log r_i \wedge d\gamma_i.$$
\end{prop}
Let $y_i=\log r_i$, and we see that $\omega^\vee=(2\pi)^n\sum d
y_i\wedge d\gamma_i$. Hence $M$ carries the symplectic structure
of the cotangent bundle over $(S^1)^n$.  The coordinates
$\gamma_i$ are on the base $(S^1)^n$ while $y_i$ live on the
fiber. The coordinate systems $(r,\gamma)$ and $(y,\gamma)$ are
used throughout the paper with the relation $y_i=\log r_i$. We
sometimes write $T^*(S^1)^n$ as a synonym for the moduli space
$M$, since we primarily consider the symplectic structure of $M$.
For the complex structure on $M$, there is
\begin{prop}
The moduli space $M$ is biholomorphic to a subset of $(\mathbb
C^*)^n$, given by the complex coordinates
$z_j=\exp(-2\pi\phi_j(L))\hol_\nabla([S^1(r_j)])$. The map
$\phi_j$ is the $j$-th component of the moment map $\phi$, while
$\hol_\nabla([S^1(r_j)])$ is the holonomy of $\nabla$ with respect
to $[S^1(r_j)]$.
\end{prop}
A straightforward calculation shows that
$$z_j(L,\nabla)=\exp(-\frac{2\pi r_j^2}{1+\sum\limits_{i=1}^n r_i^2}+2\pi i \gamma_j).$$ It can be verified that $\omega^\vee$ is indeed a K\"ahler
structure with respect to the complex structure.

The superpotential of $M$ is obtained by $m_0$ obstruction to
Floer homology, roughly speaking, counting holomorphic discs
attached to a special Lagrangian. For here, it is explicitly given
by
\begin{prop}
The superpotential $W$ on $M$ is given by
$$W=z_1+\dots+z_n+\frac{e^{-2\pi}}{z_1 z_2\dots z_n}.$$
\end{prop}

\section{T-duality and constructible sheaves}

\subsection{T-duality on torus fibers}

As the definition of $M$ shows, we can go from a \textit{flat
$U(1)$ gauge field} on a special Lagrangian fiber in $\mathbb
P^n\backslash D$ to a \textit{point} in the corresponding dual
fiber in $M$. Leung-Yau-Zaslow and Arinkin-Polishchuk apply a
similar transformation in \cite{lyz} and \cite{ap} respectively,
and go from a Lagrangian to a gauge field. Here we do this in the
other way, namely, from a gauge field to a Lagrangian. From any
holomorphic line bundle on $\mathbb P^n$, we construct an exact
Lagrangian in $M$.

We endow a line bundle $E$ on $\mathbb P^n$ with a $T^n$-invariant
hermitian metric $h$, constant on each fiber of $\phi$. The
canonical connection on $E$ with respect to this metric is
$\nabla_{E,h}=d-iA_{E,h}$, with the connection $1$-form
$A_{E,h}=i\partial h\cdot h^{-1}$ in some
trivialization.\footnote{We use the physics notation, such that
$A_{E,h}$ is real-valued.} The restriction of $\nabla_{E,h}$ to
any fiber $L$ of $\phi|_{\mathbb P^n\backslash D}$ gives rise to a
connection $\nabla_{E,h}|_L$ on the special Lagrangian $L$.
\begin{lemm}
$\nabla_{E,h}|_L$ is a flat connection over $L$.
\end{lemm}
\proof The connection $\nabla_{E,h}$ can be written as
$d-iA_{E,h}$. The connection $1$-form $A_{E,h}$ is given by
\begin{eqnarray*}
A_{E,h}&=&i\partial h\cdot h^{-1}\\
       &=&-h^{-1}\cdot\sum_{i=1}^n \frac{\partial h}{\partial
       r_i}r_id\theta_i+\text{terms in }dr_1,\dots,dr_n.
\end{eqnarray*}
The metric $h$ is a function of $r=(r_1,\dots,r_n)$ and it does
not depend on $\theta$, since it is constant on each fiber. Hence
the restriction of $A_{E,h}$ on each fiber gives vanishing
curvature.\qed

Recall from Section 1 that $M$ is the space of non-singular
$T^n$-orbits together with flat connections, and $M\cong (\mathbb
C^*)^{n}$ as a symplectic manifold with coordinates $(r,\gamma)$
and the symplectic form $$\omega^\vee=(2\pi)^n \sum_{i=1}^n
d\log(r_i)\wedge d\gamma_i.$$ We define the submanifold $\mathcal
L(E, h)\subset M$ to be
$$\{(L,\nabla_{E,h}|_{L}):\text{$L$ is an $T^n$-orbit in $\mathbb
P^n\backslash D$}\}.$$ In coordinates, we see that
$\nabla_{E,h}|_{L}$ is smooth as a function of the fiber $L$, and
hence $\mathcal L(E,h)$ is a submanifold of $M$.
\begin{prop}
$\mathcal L(E,h)$ is an exact Lagrangian submanifold.
\end{prop}
\proof In the coordinates $(r,\gamma)$ of $M$, $$\mathcal
L(E,h)=\left\{(r_1,\dots,r_n;-h^{-1}\cdot\frac{\partial h}{\partial
       r_1}r_1,\dots, -h^{-1}\cdot\frac{\partial h}{\partial
       r_n}r_n), (r_1,\dots,r_n)\in(\mathbb R^+)^n\right\}.$$
The tangent space of $\mathcal L(k,h)$ at any point is spanned by
the collection $\{\partial_{r_i}-h^{-1}\cdot\sum_{j=1}^n r_j\frac{\partial^2
h_j}{\partial r_i\partial r_j}\partial_{\gamma_j}, 1\le i\le n\}$.
We have \begin{eqnarray*}&\
&\omega^\vee(\partial_{r_i}-h^{-1}\cdot\sum_{j=1}^n r_j\frac{\partial^2
h_j}{\partial r_i\partial
r_j}\partial_{\gamma_j},\partial_{r_{i'}}-h^{-1}\cdot\sum_{j'=1}^n
r_{j'}\frac{\partial^2 h_{j'}}{\partial r_{i'}\partial
r_{j'}}\partial_{\gamma_{j'}})\\&=&(2\pi)^nh^{-1}(-\frac{\partial^2
h_j}{\partial r_i\partial r_{i'}}+\frac{\partial^2 h_j}{\partial
r_i\partial r_{i'}})\\&=&0.\end{eqnarray*} The
simple-connectedness of $\mathcal L(E,h)$ implies it is exact.\qed

\subsection{Objects in the Fukaya category}

In this subsection, we show that for the canonical metric $h_k$ on
$\mathcal O(k)$, the Lagrangian $\mathcal L(\mathcal O(k),h_k)$
can be endowed with a canonical brane structure, thus it is an
object in the Fukaya category $Fuk(M)$.

For the line bundle $\mathcal O(1)$, on the open set
$U=\{(1:z_1:\dots:z_n)\}\subset \mathbb P^n$ we can write any
$x\in \mathcal O(1)|_U$ as $\{(1:z_1:\dots:z_n),\xi\}$ by a local
trivialization of $\mathcal O(1)$. The natural $T^n$-invariant hermitian metric $h_1$ on $\mathcal O(1)$ is given by
$$h_1(x,y)=\frac{\langle\xi,\eta\rangle}{1+\sum\limits_{i=1}^{n}
|z_i|^2}.$$ 

Writing $z_i=r_ie^{i\theta_i}$, the canonical
connection $\nabla_{\mathcal O(1), h_1}=d-iA_{\mathcal O(1), h_1}$
on $\mathcal O(1)$ with respect to $h_1$ is
\begin{eqnarray*}
A_{\mathcal O(1), h_1}&=&i\partial h_1\cdot
h_1^{-1}\\&=&\frac{r_1^2d\theta_1}{1+\sum\limits_{i=1}^n r_i^2}
+\dots+\frac{r_n^2d\theta_n}{1+\sum\limits_{i=1}^n
r_i^2}+\textit{\ terms of\ }dr_1,\dots,dr_n.\end{eqnarray*} Hence
the Lagrangian
\begin{eqnarray*}&\ &\mathcal L(\mathcal O(1), h_1)\\ &=& \left\{(r,
\gamma^{(1)}(r)):r=(r_1,\dots,r_n)\in (\mathbb R^+)^n,
\gamma^{(1)}(r)=(\frac{r_1^2}{1+\sum\limits_{i=1}^n r_i^2}
,\dots,\frac{r_n^2}{1+\sum\limits_{i=1}^n r_i^2})\right\},
\end{eqnarray*}
which is obviously the graph of the $(S^1)^n$-valued function
$\gamma^{(1)}$.
 For
any other holomorphic line bundle $\mathcal O(k)$, let $h_k$
denote $(h_1)^k$, and this construction gives rise to a Lagrangian
$$\mathcal
L(\mathcal O(k), h_k)=\{(r, \gamma^{(k)}):r=(r_1,\dots,r_n)\in
(\mathbb R^+)^n, \gamma^{(k)}=k\gamma^{(1)}\}.$$

\begin{center}
\includegraphics{fig.1}
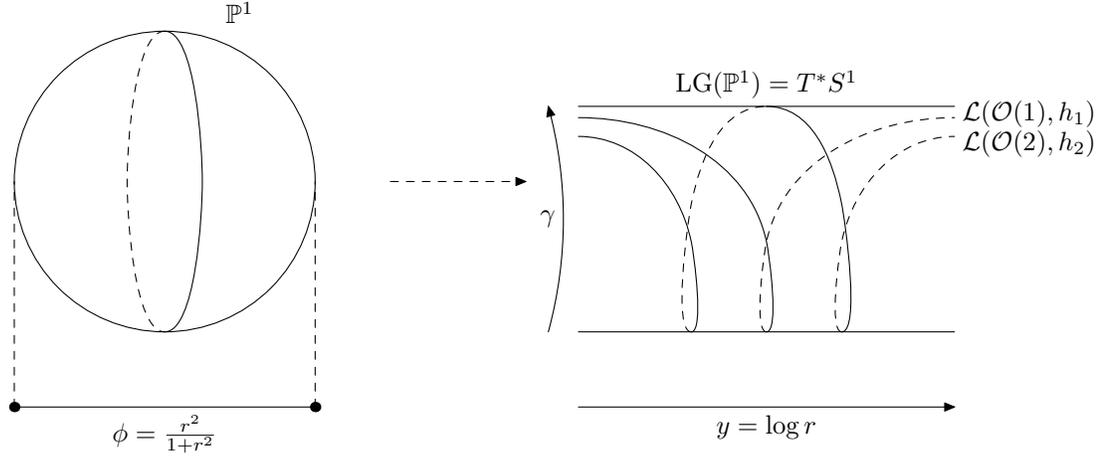
\parbox{8cm}{Fig.1 Fiberwise T-duality transformation for $\mathbb
P^1$. The Lagrangians $\mathcal L(\mathcal O(1), h_1)$ and
$\mathcal L(\mathcal O(2), h_2)$ shown are obtained from line
bundles $\mathcal O(1)$ and $\mathcal O(2)$ respectively.}
\end{center}

We adapt the definition of the Fukaya category $Fuk(M)$ for
cotangent bundles from \cite{nz}. The moduli space $M$ is
symplectomorphic to the cotangent bundle $M=T^* (S^1)^n$, where
$S^1=\mathbb R/\mathbb Z$. Moreover, it is already equipped with a
standard symplectic form $\omega^\vee=(2\pi)^n\sum d y_j\wedge
d\gamma_j$. The Fukaya category of a cotangent bundle has been defined in
\cite{nz}, \footnote{This definition differs from the ``wrapped"
category of Fukaya-Seidel-Smith \cite{fss}.} and we apply that definition to $M=T^*(S^1)^n$.
The variable $\gamma_j$ are the coordinates on the base which are periodic with period 
$1$ and $y_j$ are the coordinates on the fiber. Fix a metric on
the base $(S^1)^n$ to be
$g=d\gamma_1^2+d\gamma_2^2+\dots+d\gamma_n^2$. We would like to
employ the result of \cite{nz} to perform calculations in the Fukaya category.
The projective space $\mathbb P^n$ has an exceptional collection.
We choose one here: $\mathcal O(-n-1),\mathcal
O(-n),\dots,\mathcal O(-1)$. The philosophy is to show that the
Lagrangians $\{\mathcal L(\mathcal O(k),h_k)\}_{k=-n-1}^{-1}$
constructed from this collection form a derived Fukaya subcategory
equivalent to $D^bCoh(\mathbb P^n)$. From now on, we only consider
these objects, as well as their cones, shifts and sums in the Fukaya category of $M$. 
Let us recall some basic facts concerning the geometry of cotangent bundles from
\cite{nz}.

The next lemma shows that $\mathcal L(\mathcal O(-1), h_{-1})$ is
a graph over an open set.
\begin{lemm} Let $\mathcal T=\{(\gamma_1,\dots,\gamma_n)|\gamma_i<0\text{, and
}\sum_{i=1}^n \gamma_i>-1\}$, an $n$-cell in the base $(S^1)^n$ of
$M=T^*(S^1)^n$. The Lagrangian $\mathcal L(\mathcal
O(-1),h_{-1})$ is the graph $\Gamma_{df}$ of an exact one form $df$ for some $f:\mathcal
T\rightarrow \mathbb R$.
\label{graph}
\end{lemm}
\proof We know that $\mathcal L(\mathcal O(-1),h_{-1})$ is given
by $n$ equations
$$\gamma_j=-\frac{r_j^2}{1+\sum\limits_{i=1}^n r_i^2}$$ for $1\le
j\le n$ and $(r_1,\dots,r_n)\in (\mathbb R^+)^n$. Rewriting $r_i$
in the form of $\gamma_j$, we get
$$r_i=(-\frac{\gamma_j}{1+\sum\limits_{j=1}^n \gamma_j})^{1/2}$$
for $1\le i \le n$ and $(\gamma_1,\dots,\gamma_n)\in \mathcal T$.
Using the coordinate system $(y,\gamma)$ on the cotangent bundle,
where $y_i=\log r_i$, we find
$$
y_i=\frac{1}{2}
\log(-\frac{\gamma_i}{1+\sum\limits_{j=1}^n\gamma_j}).$$ Now note
$$y_i=\frac{\partial f}{\partial \gamma_i},$$
where
$$f=\frac{1}{2}\sum_{i=1}^n\gamma_i\log(-\gamma_i)-\frac{1}{2}(1+\sum_{j=1}^n\gamma_j)\log(1+\sum_{j=1}^n\gamma_j).$$
\qed

\begin{center}
\includegraphics{fig.2}\\
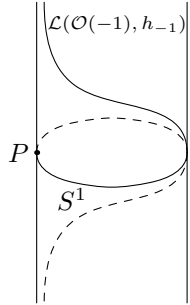
\parbox{8cm}{Fig.2 The Landau-Ginzburg mirror of $\mathbb
P^1$. The Lagrangian $\mathcal L(\mathcal O(-1), h_{-1})$ is a
graph over the open interval $\mathcal T=S^1\backslash P$.}
\end{center}

\begin{lemm} There are canonical brane structures for
objects $\mathcal L(\mathcal O(k),h_k)$ for $k\in
\{-n-1,\dots,-1\}$. Hence we have obtained $n$ objects $\mathcal
L(\mathcal O(k),h_k)$ in the Fukaya category $Fuk(M)$.
\end{lemm}
\proof The Lagrangian $\mathcal L(\mathcal O(-1),h_{-1})$ is
canonically Hamiltonian isotopic to $\mathcal T$, inside
$T^* (S^1)^n|_{\mathcal T}$ (this is a hamiltonian isotopy inside $T^* (S^1)^n|_{\mathcal T}$ but not the whole $T^*(S^1)^n$). Let $p: T^*(S^1)^n\rightarrow
(S^1)^n$ be the projection to the base. This Hamiltonian isotopy
can be achieved by the Hamiltonian flow $\varphi_{H,t}$, where $H=
f \circ p$, which takes $\mathcal L(\mathcal O(-1),h_{-1})$ to
$(1-t)\mathcal L(\mathcal O(-1),h_{-1})$. In particular, when
$t=1$, one arrives at $\mathcal T$. We can equip $\mathcal T$ with
grading $0$. Because $\mathcal T$ is canonically Hamiltonian
isotopic to $\mathcal L(\mathcal O(-1),h_{-1})$, there is a
canonical grading for $\mathcal L(\mathcal O(-1),h_{-1})$. As for
the pin structure, since $\mathcal L(\mathcal O(-1),h_{-1})$ is
contractible it has a trivial pin structure. So we obtain a
canonical brane structure on $\mathcal L(\mathcal O(-1),h_{-1})$.

Let $\widetilde S^1=\mathbb R/(n+1)\mathbb Z$, and let $\widetilde
M=T^*(\widetilde S^1)^n$. Consider the $(n+1)^n$-covering $\pi:
\widetilde M\rightarrow M$, given by $ \pi: (y, \tilde\gamma \mod
(n+1))\mapsto (y,\tilde\gamma \mod 1).$ The variables
$\tilde\gamma=(\tilde\gamma_1,\dots,\tilde \gamma_n)$ are $n+1$
periodic in each component. To treat $\mathcal L(\mathcal
O(k),h_k)$ for $k\le -2$, we need to consider the lifts of these
Lagrangians under this covering map $\pi$. These Lagrangians
become graphs over open sets in $\widetilde M$ after the lifting.

Let $a=(a_1,\dots,a_n)\in (\mathbb Z/(n+1))^n$, and assume each
$a_i$ takes integer value from $-n$ to $0$. We can define open
sets $$\mathcal
U(k)^{(a_1,\dots,a_i)}=\left\{\tilde\gamma|\tilde\gamma_i<a_i,
\sum_{j=1}^n \tilde\gamma_j > k + \sum_{j=1}^n a_j\right\}$$ for
$k\in\{-n-1,\dots,-1\}$. An exact Lagrangian $\mathcal L(\mathcal
O(k),h_k)$ admits $(n+1)^n$ possible lifts. Let $f^0_{-1}$ be a
real function on $\mathcal U^0(-1)$, defined by:
$$f_{-1}^0=\frac{1}{2}\sum_{i=1}^n\tilde\gamma_i\log(-\tilde\gamma_i)-\frac{1}{2}(1+\sum_{j=1}^n\tilde\gamma_j)\log(1+\sum_{j=1}^n\tilde\gamma_j).$$
Note that $f^0_{-1}$ is one of the lifts of $f:\mathcal
T\rightarrow \mathbb R$. For $k\in\{-n-1,\dots,-1\}$, there are
$(n+1)^n$ lifts of $\mathcal L(\mathcal O(k),h_k)$, namely, given
by the graph of the differential of
$$f^a_{k}=f^0_{-1}\left(\frac{\tilde\gamma_1-a_1}{-k},\frac{\tilde\gamma_2-a_2}{-k},\dots,\frac{\tilde\gamma_n-a_n}{-k}\right),$$
for any $a=(a_1,\dots,a_n)\in (\mathbb Z/(n+1))^n$. We denote the
lifted Lagrangian as the differential of the above function by
$\mathcal L^a(\mathcal O(k),h_k)$.

\begin{center}
\includegraphics{fig.3}
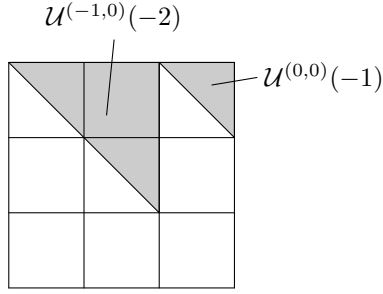
\parbox{8cm}{Fig.3 This square represents the base $(\widetilde S^1)^2$ of $T^*(\widetilde S^1)^2$, as the lift $\widetilde M$ of the Landau-Ginzburg mirror $M$ of
$\mathbb P^2$. The Lagrangian $\mathcal L^{(0,0)}(\mathcal
O(-1),h_{-1})$ is a graph over the small shaded triangle $\mathcal
U^{(0,0)}(-1)$, while $\mathcal L^{(-1,0)}(\mathcal O(-2),h_{-2})$
is a graph over $\mathcal U^{(-1,0)}(-2)$.}
\end{center}

Since the topology of $\mathcal L(\mathcal O(k),h_k)$ is trivial,
it has a trivial (and canonical) pin structure. We know that
naturally $\mathcal L^0(\mathcal O(k),h_k)$ has a canonical brane
structure, by the same argument for $\mathcal L(\mathcal
O(-1),h_{-1})$. The covering map $\pi$ acts trivially on the
phase, and hence we can make $\pi$ into a \textit{graded} covering
$\widetilde \pi$ with trivial grading. Under the graded covering map
$\widetilde \pi$, the natural grading of $\mathcal L^0(\mathcal
O(k),h_k)$ is mapped to a grading of $\mathcal L(\mathcal
O(k),h_k)$ in $M$, giving a canonical brane structure for
$\mathcal L(\mathcal O(k),h_k)$. Notice that our construction does
not depend on the lift of $\mathcal L(\mathcal O(k),h_k)$. If
$\mathcal L^a(\mathcal O(k),h_k)$ is another lift of $\mathcal
L(\mathcal O(k),h_k)$, for any $x \in \mathcal L^0(\mathcal
O(k),h_k)$ and $x' \in \mathcal L^a(\mathcal O(k),h_k)$ such that
$\pi (x)= \pi (x')$, we have $\widetilde \alpha_{\mathcal L^0(\mathcal
O(k),h_k)}(x)=\widetilde \alpha_{\mathcal L^a(\mathcal
O(k),h_k)}(x')$, where $\widetilde \alpha_{\mathcal L^0(k,h_k)}$ and
$\widetilde \alpha_{\mathcal L^a(\mathcal O(k),h_k)}$ are canonical
gradings of $ {\mathcal L}^0(\mathcal O(k),h_k)$ and ${\mathcal
L}^a(\mathcal O(k),h_k)$ respectively. Hence different lifts give
the same grading for $\mathcal L(\mathcal O(k), h_k)$. \qed
\begin{rema}
Although we are only worrying about finitely many $\mathcal
L(\mathcal O(k),h_k)$ for $k\in\{-n-1,\dots,-1\}$, this lemma
actually holds for all $k$, i.e. any $\mathcal L(\mathcal
O(k),h_k)$ for $k\in \mathbb Z$ has a canonical brane structure
and can be made into an object in the Fukaya category $Fuk(M)$. Therefore from any line bundle $\cO(k)$ on $\bP^n$, we can construct a Lagrangian brane $\cL(\cO(k), h_k)$ on the mirror side. Moreover, this construction does not essentially depend on the choice of $T$-invariant metric, although we are using the canonical metric $h_k$ here. Different metrics give rise to quasi-isomorphic branes in the Fukaya category. Let $h'_k=e^\lambda h_k$ be another metric on $\cO(k)$ where $\lambda$ is a $T^n$-invariant function on $\bP^n$. A straightforward calculation shows that the Lagrangian $\cL(\cO(k), h'_k)=\phi_1(\cL(\cO(k), h_k)$, where $\phi$ is the hamiltonian flow generated by the function $\lambda$. A more detailed argument of non-characteristic isotopy in \cite{n} shows that $\cL(\cO(k), h'_k)$ and $\cL(\cO(k), h_k)$ are quasi-isomorphic in the Fukaya category $Fuk(M)$. For a direct treatment of an arbitrary $T$-invariant metric in the T-duality, please see \cite{fltz}. In the rest of this paper, we simply denote $\cL(\cO(k), h_k)$ by  $\cL(k)$ for convenience.
\end{rema}

When considering homological mirror symmetry, one actually deals
with the derived version of the triangulated envelope of the
Fukaya category. There are several ways to define the triangulated
envelope of a Fukaya category. Here we adopt the method of Yoneda
embedding, which agrees with the definition in \cite{n}.

For any $A_\infty$ category $\cA$, the Yoneda embedding $\mathcal Y: \cA \rightarrow
mod(\cA)$ maps an object $L\in \cA$ to an
$A_\infty$-module $hom_{\cA}(-,L)$. We write $Tr(\cA)$ for the category of twisted complexes of modules in $\cY(\cA)$ as a version of the triangulated envelope of $\cA$.

We study the Fukaya category $\mathcal F$ as a full sub-category
of $Fuk(M)$ containing objects $\mathcal L(k)$, $-n-1\le k
\le -1$ in this paper. Precisely,
\begin{defi}
$\mathcal F$ is the full $A_\infty$ subcategory of $Fuk(M)$ consisting
of $n$ objects $\mathcal L(k)$ where $k\in\{-n-1,\dots,-1\}$.
\end{defi}
The derived category of ${\mathcal F}$, denoted by $D \mathcal F$,
is a triangulated category $H^0(Tr({\mathcal F}))$. Note that $D
\mathcal F\subset D Fuk(M)$ is a full subcategory. We state our
main theorem:
\begin{theo*}
The (bounded) derived Fukaya category $D {\mathcal F}$ is equivalent to $D^b
Coh(\mathbb P^n)$, the derived category of coherent sheaves on
$\mathbb P^n$.
\end{theo*}

\subsection{Passing to standard branes}

Let $X$ be a real analytic manifold. A standard brane in the Fukaya category $Fuk(T^*X)$ over an open set $U\subset X$ is the graph of the differential $d\log m$ equipped with the canonical brane structure, where $m$ is a defining function of $\partial U$ on $\overline U$:  $m=0$ on $\partial U$ and $m>0$ on $U$. The quasi-isomorphism class of this brane does not depend on the choice of the particular $m$. 

We would like to remind that the Lagrangian $\cL(-1)$ is the graph of the exact differential form $df$ on the open set $\cT\subset (S^1)^n$. Note that we cannot claim $\mathcal L(-1)$ is a standard brane in $Fuk(M)$ at present, since $e^f$ does not go to $0$ near the boundary. However, it looks very much like a standard Lagrangian, i.e. the covector $df$ points inward near the boundary $\partial \mathcal T$ and its length is arbitrarily large. This section shows that $\mathcal L(-1)$ is indeed isomorphic to a standard brane over the open set $\cT$, which allows us to apply the microlocalization functor in \cite{nz} to pass into the category of constructible sheaves.

For any real analytic manifold $X$, let $Sh_{naive}(X)$ be the
triangulated dg category whose objects are complexes of sheaves
with bounded constructible cohomology, and whose morphisms are the
usual complexes of morphisms. Then we take $Sh(X)$ be the dg
quotient of $Sh_{naive}(X)$ with respect to the subcategory
$\mathcal N$ of acyclic objects \cite{ke}. 

The result of \cite{nz} says that there is an $A_\infty$ microlocalization functor from $Sh(X)$ to
$TrFuk(T^*X)$,
such that the induced functor for derived categories $ D
Sh(X)\rightarrow DFuk(T^*X)$ is an embedding. The
functor is generated by sending $i_* \mathbb C_U$ for any
open set $i: U\hookrightarrow X$ to the standard Lagrangian over $U$. Particularly, for $X=(S^1)^n$, we use $\mu$ as the microlocalization functor from $Sh((S^1)^n)\rightarrow TrFuk(M)$ or the derived version $DSh((S^1)^n)\rightarrow DFuk(M)$, depending on the context. Similarly, there is a microlocaliztion functor $\widetilde \mu:Sh((\widetilde S^1)^n)\rightarrow TrFuk(\widetilde M)$. Due to this fact, we denote the standard brane over $U$ in $Fuk(M)$ by $\mu(i_*\bC_U)$, and the standard brane over $\widetilde U$ in $Fuk(\widetilde M)$ by $\widetilde \mu(i_*\bC_{\widetilde U})$. 

The \emph{normalized geodesic flow} $\varphi_t$ perturbs the objects when defining the morphisms in $Fuk(M)$. Since we've already chosen a metric on the base $(S^1)^n$ of $M=T^*(S^1)^n$, i.e $g=d\gamma_1^2+\dots+d \gamma_n^2$, the normalized geodesic flow $\varphi$ is
$$\varphi_t(y,\gamma)= (y, \gamma+t \frac{y^*}{\|y\|}),$$
where $y^*\in T_\gamma (S^1)^n$ is the dual of $y\in T^*_\gamma (S^1)^n$ with respect to this metric. Note that this flow is only defined on $(T^*(S^1)^n)_0=\{(y,\gamma)\in M| y\ne 0\}\subset M$, i.e. away from the zero section in the cotangent bundle. Let $s$ be an arbitrary point in $\partial \cT$, and $L_{\{s\}*}$ be the brane supported on the fiber Lagrangian $T^*_s (S^1)^n$. Define
\begin{eqnarray*}
\cL(-1)_0= \cL(-1)\cap (T^*(S^1)^n)_0,\\
(L_{\{s\}*})_0=L_{\{s\}*}\cap (T^*(S^1)^n)_0.
\end{eqnarray*}

\begin{lemm}
There exists a $\delta>0$, such that
$$0\le t_1\le t_2<\delta\Rightarrow \varphi_{t_1} ((L_{\{s\}*})_0)\cap \varphi_{t_2}(\cL(-1)_0)=\emptyset.$$
\label{nohom}
\end{lemm}
\proof
The Lagrangian $\cL(-1)$ is the graph of $df$ over the open set $\cT\subset (S^1)^n$. The function $f$ here, as given in Lemma \ref{graph}, is
$$f=\frac{1}{2} \sum\limits_{i=1}^n \gamma_i \log (-\gamma_i)-\frac{1}{2}(1+\sum\limits_{j=1}^n \gamma_j)\log(1+\sum\limits_{j=1}^n \gamma_j).$$
Here we assume $\gamma_i$ takes value in $(-1,0)$. Denote $Z_i=\gamma_i$ for $1\le i \le n$, and $$Z_0=-1-\sum\limits_{j=1}^n \gamma_j.$$Hence the function
$$f=\frac{1}{2}\sum_{i=0}^n Z_i\log(-Z_i),$$
while the open set $\cT$ is characterized by $$\cT=\{Z_i<0|0\le i \le n\}.$$

Since $s\in \partial \cT$, there exists a non-empty subset of index $I\subset \{0,\dots, n\}$, such that $Z_i(s)=0$ for $i\in I$ and $Z_i(s)< 0$ for $i\in \{0,\dots, n\}\backslash I$. Note that $1\le |I|\le n$ because of the constraint $\sum_{i=0}^n Z_i=-1$ so that we cannot have $Z_i(s)=0$ for all $0\le i\le n$. There is a bound $R>0$ such that $|\log(-Z_i(s))|<R/2$ for all $i\in \{0,\dots,n\}\backslash I$. Choose a subset $K\subset \{0,\dots,n\}$ such that $|K|=n$ and $I\subset K$. Therefore $\{Z_i|i\in K\}$ form a coordinate system on $(S^1)^n$. Let $i_0=\{0,\dots, n\}\backslash K$ be the single index that is not in $K$. It is easy to see that the metric $$g_K=\sum_{i\in K} dZ_i^2$$ is equivalent to the standard metric $$g=d\gamma_1^2+\dots+d\gamma_n^2=\sum_{i=1}^n dZ_i^2.$$ Therefore there is a bound $Q$ with
$$\|\cdot\|_{g_K}\le 1/Q\|\cdot\|_g.$$

For any $M>0$, there is a $\delta>0$, such that for any $\gamma\in \cT$ with $\dist_g(\gamma,s)<\delta$, $\log (-Z_i(\gamma))<-M$ for all $i\in I$ and $|\log(-Z_i(\gamma))|<R$ for $i\in \{0,\dots, n\}\backslash I$. We choose a large $M$ such that $\frac{2nR}{Q^2(M-R)}<1$, and a $\delta$ corresponding to this $M$.

For any $(y,\gamma)\in \cL(-1)$ and any given $0\le t_1\le t_2<\delta$, it suffices to show that $\varphi_{t_1}(y, s)\ne \varphi_{t_2}(y,\gamma)$ to finish the proof. We know that
$$\varphi_{t_1}(y, s)=(y, s+t_1\frac{y^*}{\|y\|_g}),\text{and }\varphi_{t_2}(y, \gamma)=(y,\gamma+t_2\frac{y^*}{\|y\|_g}).$$
Therefore $\varphi_{t_1}(y, s)=\varphi_{t_2}(y,\gamma)$ implies that $v=(t_2-t_1)\frac{y^*}{\|y\|_g}$, where $v=s-\gamma$ is considered as a vector in $T(S^1)^n$. Decompose $v=\sum_{i\in K} v_i\partial_{Z_i}$. When $i\in I$, $$
Z_i(s)=0, Z_i(\gamma)<0\Rightarrow v_i>0.$$
We will show that $v=(t_2-t_1)\frac{y^*}{\|y\|_g}$ is impossible.

\smallskip

\noindent
Case 1.  $\dist_g(\gamma, s)<\delta$. We have 
\begin{eqnarray*}
y=df(\gamma)&=&\sum_{i\in K}\frac{\partial f}{\partial Z_i} dZ_i\\
			&=&\sum_{i\in K} \frac{1}{2}(\log(-Z_i)-\log(-Z_{i_0}))dZ_i.
\end{eqnarray*}
For $i\in I$, 
$$\frac{\partial f}{\partial Z_i}<\frac{1}{2}(-M+R)<0,$$
while for $i\in K\backslash I$,
$$|\frac{\partial f}{\partial Z_i}|<\frac{1}{2}R<R.$$
The length of $y$ satisfies
$$
\|y\|_g\ge Q\|y\|_{g_K}\ge Q\sqrt{\sum_{i\in I} \left | \frac{\partial f}{\partial Z_i}\right |^2}\ge Q\sqrt{\frac{1}{4}(M-R)^2}\ge \frac{1}{2}Q(M-R).
$$

Therefore
\begin{eqnarray*}
\langle v, \frac{y}{\|y\|_g} \rangle &<&\sum_{i\in I} v_i \frac{-M+R}{Q(M-R)}+\sum_{i\in K\backslash I} |v_i|\frac{ 2R}{Q(M-R)} \\
								     &<& \frac{2R}{Q(M-R)}\|v\|_{g_K} |K\backslash I|< \frac{2nR\|v\|_g}{Q^2(M-R)}.
\end{eqnarray*}

Since $\frac{2nR}{Q^2(M-R)}<1$, the inner product
$$\langle v, \frac{y^*}{\|y\|_g} \rangle_g=\langle v, \frac{y}{\|y\|_g} \rangle <\|v\|_g.$$
Note that the dual $y^*$ is taken with respect to the standard metric $g$. We know the length of $\frac{y}{\|y\|_g}$ is $1$, and this shows that $y^*$ is not parellel and in the same direction with $v$. Thus the equality $v=(t_2-t_1)\frac{y^*}{\|y\|_g}$ is impossible.

\smallskip

\noindent
Case 2. $\dist_g(\gamma,s)\ge\delta$. Therefore  $v=(t_2-t_1)\frac{y^*}{\|y\|_g}$ is impossible to hold since $\|v\|_g\ge\delta$ while the length of the right hand side is $t_2-t_1<\delta$.

\qed

\begin{lemm}
The Lagrangian brane $\mathcal L(-1)$ is isomorphic to a standard
Lagrangian over the open set $\mathcal T$ in the category
$DFuk(M)$. Similarly, all $\mathcal L^a(k)$ are
isomorphic to standard Lagrangians over $\mathcal U^a(k)\subset
(\widetilde S^1)^n$ in $DFuk(\widetilde M)$.\label{standard}
\end{lemm}
\proof Since the functions $f$ and $f^a_k$ are essentially the
same, i.e. only differing by scaling of the domain, it suffices to
show this lemma only for $\mathcal L(-1)$.

To prove that $\cY(\cL(-1))\cong \cY(\mu(i_*\bC_\cT))$, 
we first fix a triangulation $\Lambda$ of the base $(S^1)^n$
containing $\cT$ and each stratum of its boundary.
The technique of \cite{n} exploits the triangulation to resolve the
diagonal standard, i.e. the identity functor.
What emerges is that the Yoneda module of any object $\mathcal Y(L)$ is expressed
in terms of (sums and cones of shifts of) Yoneda modules from standards,
$\mathcal Y(\mu(i_*\bC_{T})),$ where $T\in \Lambda.$
The coefficient of the Yoneda standard module $\mathcal Y(\mu(i_*\bC_{T})),$
takes the form $hom_{DFuk(M)}(L_{\{s\}*},L),$
where $s$ is any point in $T$ (contractibility of $T$ means
that the choice is irrelevant up to isomorphism) -- see Proposition 4.4.1 and Remark 4.5.1 of \cite{n}.

Now apply this to $\cL(-1).$ Note that $\Lambda$ contains all strata of $\cT$ and $\partial \cT$.
Let $T\neq \cT$ and let $s\in T.$ 
Then if $\overline{T}\cap \cT = \emptyset,$
clearly $hom_{DFuk(M)}(L_{\{s\}*},\cL(-1)) = 0,$
since $L_{\{s\}*}$ is just the fiber $T^*_s(S^1)^n.$
Otherwise, if
$T\cap \partial \cT$ is nonempty, then Lemma \ref{nohom}
ensures us that the hom space $hom_{DFuk(M)}(L_{\{s\}*},\cL(-1)) = 0.$
Finally, if $T = \cT,$ then since $\cL(-1)$ is
a graph over $T,$ the morphism space
$hom_{DFuk(M)}(L_{\{s\}*},\cL(-1)) = \bC.$
Therefore, $\cY(\cL(-1)) \cong \cY(\mu(i_*\bC_\cT)).$
Note that the result is independent of how $\Lambda$ was chosen.\qed

\section{The quasi-equivalence of the categories}

Recall $\widetilde M=T^*(\widetilde S^1)^n$, where $\widetilde
S^1$ is identified with $\mathbb R/(n+1)\mathbb Z$. The variable
$\tilde\gamma=(\tilde\gamma_1,\tilde\gamma_2,\dots,\tilde\gamma_n)$
on the base is defined mod $n+1$, and $y\in
T^*_{\tilde\gamma}(\widetilde S^1)^n\cong\mathbb R^n$ is the
variable in the fiber. The covering map $\pi:\widetilde
M\rightarrow M$ is given by $(y, \tilde\gamma\mod (n+1))\mapsto
(y,\tilde\gamma\mod 1)$.

The Fukaya category ${\mathcal F}$ admits a lift $\widetilde
{\mathcal F}$ through the covering map $\pi$. The objects of
$\widetilde {\mathcal F}$ are all lifts of each object in
${\mathcal F}$. A lift $\mathcal L^a(i)$ of $\mathcal L(i)$ is the
graph of a differential $1$-form over an open set $$\mathcal
U(i)^a=\{\tilde\gamma|\tilde\gamma_j<a_j, \sum_j\tilde\gamma_j > i
+ \sum_j a_j\},$$ where $a=(a_1,\dots,a_n)$, $a_j\in \mathbb
Z/(n+1)$. Here we assume $a_j$ takes an integer value ranging from
$-n$ to $0$. The morphisms of $\widetilde {\mathcal F}$ are
inherited from $Fuk(\widetilde M)$ as
$$hom_{\widetilde{\mathcal F}}(\mathcal L^a (i), \mathcal L^b
(j))= hom_{Fuk(\widetilde M)} (\mathcal L^a
(i),\mathcal L^b(j)).$$ The composition maps $m_k$ are the same as
in $Fuk(\widetilde M)$. The category
${\widetilde{\mathcal F}}$ is a full $A_\infty$-subcategory of
$Fuk(\widetilde M)$.

Let $\widetilde{\mathcal D}$ be the differential graded category
containing objects $i_* \mathbb C_{\mathcal U^a(i)}$ for $-n-1\le
i \le -1$ and $a=(a_1,\dots,a_n)\in (\mathbb Z/(n+1))^n$. $\widetilde \cD$ is a
full subcategory of $Sh((\widetilde S^1)^n)$. 

There is a natural $\Gamma=(\mathbb Z/(n+1))^n$ action on
$\widetilde M$ given by the deck transformation, i.e. $\alpha\in
\Gamma: (y,\tilde \gamma)\mapsto (y,\tilde \gamma+\alpha)$. This
action gives rise to actions on $\mathcal {\widetilde F}$ and on
$\widetilde {\mathcal D}$. For any $\alpha\in \Gamma$, we have for
objects
$$ \alpha \cdot {\mathcal L}^a(i)=\mathcal L^{a+\alpha}(i),\qquad \alpha\cdot i_* \mathbb C_{\mathcal
U^a(i)}=i_* \mathbb C_{\mathcal U^{a+\alpha}(i)}.$$ For morphisms,
these actions induce natural isomorphisms on morphism spaces.
\begin{eqnarray*}
\alpha: hom(\mathcal L^a(i),\mathcal L^b(j))&\rightarrow&
hom(\mathcal
L^{a+\alpha}(i),\mathcal L^{b+\alpha}(j)), \\
\alpha: hom(i_* \mathbb C_{\mathcal U^a(i)},i_* \mathbb
C_{\mathcal U^b(j)})&\rightarrow& hom(i_* \mathbb C_{\mathcal
U^{a+\alpha}(i)},i_*\mathbb C_{\mathcal U^{b+\alpha}(j)}).
\end{eqnarray*}
The action respects the compositions of corresponding morphisms.

\begin{defi}
The category $\mathcal {\widetilde F}/\Gamma$ is the quotient of
$\mathcal {\widetilde F}$ with respect to the action of $\Gamma$. It
consists of $n+1$ objects, denoted as formal orbits $(\oplus_a \mathcal
L^a(i))/\Gamma$. The morphisms are defined to be
$$hom_{\mathcal {\widetilde F}/\Gamma} ((\oplus_a \mathcal L^a(i))/\Gamma, (\oplus_b \mathcal L^b(j))/\Gamma)=(\bigoplus_{a,b}hom_{\mathcal {\widetilde F}} (\mathcal L^a(i),\mathcal
L^b(j)))/\Gamma.$$ The compositions are inherited naturally from
$\mathcal {\widetilde F}$. The quotient category $\mathcal
{\widetilde D}/\Gamma$ is defined the same way.
\end{defi}

By Lemma \ref{standard} the branes $\mathcal L^a(i)$ are quasi-isomorphic
to standards over $\mathcal U^a(i)$. Hence we have
\begin{prop} The $A_\infty$ functor $\tilde\mu$ restricted on $\cD$ gives rise to a quasi-equivalence $\tilde\mu: {\widetilde {\mathcal D}}\rightarrow \mathcal {\widetilde{\mathcal F}}$, where $\widetilde \mu$ sends $i_*
\mathbb C_{\mathcal U^a(i)}$ to $\mathcal L^a(i)$. This functor is
a quasi-isomorphism.\label{standard-like}
\end{prop}
\proof This is essentially the result of \cite{nz}, and it is
obvious that $\widetilde \mu$ is $\Gamma$-equivariant. \qed

Taking the quotient of $\mathcal {\widetilde F}$ by the action
$\Gamma$, we get an $A_\infty$ Fukaya-type category$\mathcal
{\widetilde F}/\Gamma$. This category is isomorphic to $\mathcal {
F}$.
\begin{lemm}
$\mathcal{\widetilde F}/\Gamma\cong\mathcal { F}$.
\end{lemm}
\proof Define a functor $J$ sending $(\oplus_a \mathcal
L^a(i))/\Gamma$ to $\mathcal L(i)$. Any morphism from $(\oplus_a
\mathcal L^a(i))/\Gamma$ to $(\oplus_b \mathcal L^b(j))/\Gamma$
corresponds to $(n+1)^n$ intersection points (after perturbation)
between $\mathcal L^{a_0+c}(i)$ and $L^{b_0+c}(j)$ for some $a_0$,
$b_0$ and all $c\in \Gamma$. These points are of the same degree
$d$. Under the map $\pi$, all these intersection points go to one
intersection point between $\mathcal L(i)$ and $\mathcal L(j)$ of
the same degree $d$. On the other hand, any morphism between
$\mathcal L(i)$ and $\mathcal L(j)$ can be lifted to $(n+1)^n$
intersection points between $\mathcal L^{a_0+c}(i)$ and $\mathcal
L^{b_0+c}(j)$ for some $a_0$, $b_0$ and all $c\in \Gamma$. Hence
the morphism spaces of $\mathcal {\widetilde F}/\Gamma$ and
$\mathcal { F}$ are identical.

For compositions, let us consider a polygon $Q$ bounded by
$\mathcal L(i_1),\dots,\mathcal L(i_k)$. This polygon contributes
to the composition of morphisms represented by each vertex. Since
$Q$ is a simply-connected polygon, it has $(n+1)^n$ lifts to
$\widetilde M$, the polygons $Q^c$ bounded by $\mathcal
L^{a_1+c}(i_1),\dots,\mathcal L^{a_k+c}(i_k)$, for some
$a_1,\dots,a_k$ and  all $c\in\Gamma$. Similarly, the compositions
in $\mathcal{\widetilde F}/\Gamma$ come from counting of the orbit
of polygons bounded by $\mathcal L^{a_1+c}(i_1),\dots,\mathcal
L^{a_k+c}(i_k)$ for some $a_1,\dots,a_k$ and all $c\in\Gamma$.
Hence we have proved the composition maps in both $\mathcal
{\widetilde F}/\Gamma$ and $\mathcal { F}$ are the same.\qed

Since the functor $\widetilde \mu$ is $\Gamma$-equivariant, the induced
functor $\mathcal {\widetilde D}/\Gamma\rightarrow \mathcal
{\widetilde F}/\Gamma\cong\mathcal{ F}$ is a quasi-equivalence.
Hence the functor (still denoted by $\widetilde \mu$) between the derived
categories $\widetilde \mu: D(\mathcal {\widetilde{D}}/\Gamma)\rightarrow
D\mathcal{ F}$ is a equivalence.
\begin{lemm}
$D(\mathcal {\widetilde D}/\Gamma)$ is equivalent to
$D\mathcal { F}$.\label{passtosheaves}\qed
\end{lemm}

The quiver structure of $D(\mathcal {\widetilde D}/\Gamma)$ can
be computed in a combinatorial way. Note that $D(\mathcal
{\widetilde  D}/\Gamma)$ is generated by objects in ${\mathcal
{\widetilde D}}/\Gamma$, as a triangulated category. Hence we only
deal with objects in $\mathcal {\widetilde D}$. The following
statement is the Lemma 4.4.1. from \cite{nz} concerning the
morphisms of sheaves.
\begin{lemm}
Let $X$ be a real analytic manifold. For any two open sets $i_0:
U_0\hookrightarrow X$, $i_1: U_1\hookrightarrow X$, we have a
canonical quasi-isomorphism in the dg category:
$$hom_{Sh(X)}(i_{0*}\mathbb C_{U_0},i_{1*}\mathbb C_{U_1})\simeq
(\Omega(\overline U_0\cap U_1,\partial U_0\cap U_1),d).$$ The
composition of morphisms coincides with the wedge product of
differential forms.\label{sheaftoform}
\end{lemm}

By this Lemma,
$$hom_{\widetilde {\mathcal D}}(i_*\mathbb C_{\mathcal U^a(i)}, i_* \mathbb C_{\mathcal
U^b(j)})\simeq \Omega((\overline{\mathcal U^a(i)} \cap \mathcal
U^b(j),\partial \mathcal U^a(i)\cap \mathcal U^b(j)),d).$$ Hence
$$hom_{{\widetilde{\mathcal D}}/\Gamma}(\mathcal U(i),\mathcal U(j))=(\bigoplus_{a,b}\Omega((\overline{\mathcal U^a(i)} \cap \mathcal
U^b(j),\partial \mathcal U^a(i)\cap \mathcal U^b(j)),d))/\Gamma,$$
with $\Gamma$ acting on the space of differential forms in the
obvious way. Notice that when $\mathcal U^a(i)\supset \mathcal
U^b(j)$, we have $\Omega((\overline{\mathcal U^a(i)} \cap \mathcal
U^b(j),\partial \mathcal U^a(i)\cap \mathcal
U^b(j)),d)=\Omega(\mathcal U^b(j),d)$. The cohomology
$H^*(\Omega(\mathcal U^b(j),d))=\mathbb C [0]$ ($\mathbb C$ in the
zeroeth degree). Otherwise when $\mathcal U^a(i)\not\supset
\mathcal U^b(j)$, $$\Omega((\overline{\mathcal U^a(i)} \cap
\mathcal U^b(j),\partial \mathcal U^a(i)\cap \mathcal
U^b(j)),d)\cong\Omega((\overline{\mathcal U^a(i)} \cap \mathcal
U^b(j)/(\partial \mathcal U^a(i)\cap \mathcal U^b(j)),pt),d).$$ It
follows that
\begin{eqnarray*}
&\ &H^*(\Omega((\overline{\mathcal U^a(i)} \cap \mathcal
U^b(j),\partial \mathcal U^a(i)\cap \mathcal
U^b(j)),d))\\&=&H^*(\overline{\mathcal U^a(i)} \cap \mathcal
U^b(j)/(\partial \mathcal U^a(i)\cap \mathcal U^b(j)),pt)\\&=&0.
\end{eqnarray*}
This is because when $\mathcal U^a(i)\not\supset \mathcal U^b(j)$,
$\overline{\mathcal U^a(i)} \cap \mathcal U^b(j)/(\partial
\mathcal U^a(i)\cap \mathcal U^b(j))$ is a contractible space, and
the homology is zero. Therefore, in the derived category
\begin{eqnarray*} hom_{D({\widetilde{\mathcal D}/\Gamma})}(\mathcal U(i), \mathcal
U(j)[k])&\cong&H^k((\bigoplus_{a,b}\Omega((\overline{\mathcal
U^a(i)} \cap \mathcal U^b(j),\partial \mathcal U^a(i)\cap \mathcal
U^b(j)),d))/\Gamma)\\
&=&(\bigoplus_{a,b}H^k(\Omega((\overline{\mathcal U^a(i)} \cap
\mathcal U^b(j),\partial \mathcal U^a(i)\cap \mathcal
U^b(j)),d)))/\Gamma\\
&=&\bigoplus_{\mathcal U^0(i)\supset \mathcal
U^b(j)}H^k(\Omega({\mathcal U^b(j)},d))\\
&\cong&\left\{ \begin{array}{ll} \mathbb C^{N_n(i-j)},& k=0,\\
0,& k\ne 0. \end{array}\right.
\end{eqnarray*}
Here $N_n(i-j)=\#\{b=(b_1,\dots,b_n)|b_i\le 0, \sum b_i\ge i-j
\}$, counting all possible $\mathcal U^b(j)$ in $\mathcal U^0(i)$.
Notice this explicit calculation of morphisms implies that
$\{\mathcal U(-n-1),\dots,\mathcal U(-1)\}$ is a full strong
exceptional collection of $D({\widetilde{\mathcal D}}/\Gamma)$,
since $N_n(0)=1$ and $N_n(m)=0$ when $m<0$. Let $e^b_{i,j}$ denote
the identity in $H^0(\Omega({\mathcal U^b(j)},d)$ as a
subspace of $hom_{D({\widetilde{\mathcal D}/\Gamma})}(\mathcal
U(i), \mathcal U(j))$.

We compare this exceptional collection with the exceptional
collection $\{\mathcal O(-n-1),\dots,\mathcal O(-1)\}$ of
$D^bCoh(\mathbb P^n)$. The morphism spaces of
$D(\mathcal{\widetilde D}/\Gamma)$ are
$$hom_{D({\widetilde{\mathcal D}}/\Gamma)}(\mathcal
U(i), \mathcal U(j))=\bigoplus_{\mathcal U^0(i)\supset \mathcal
U^{b}(j)}H^0(\Omega({\mathcal U^{b}(j)},d)).$$ The
constraint on the multi-index $b$ is that $$b=(b_1,\dots
b_n)\in\{b=(b_1,\dots,b_n)|b_i\le 0, \sum b_i\ge i-j \}.$$ We
construct an isomorphism $\nu$ of morphisms spaces
$hom_{D({\widetilde{\mathcal D}}/\Gamma)}(\mathcal U(i),
\mathcal U(j))$ to $hom_{D^bCoh(\mathbb P^n)} (\mathcal O(i),
\mathcal O(j))$ by
$$\begin{array}{cccc}
\nu: & hom_{D(\widetilde{\mathcal D}/\Gamma)}( \mathcal
U(i),\mathcal U(j))&\cong &hom_{D^bCoh(\mathbb P^n)}(\mathcal
O(i),\mathcal O(j))\\
\quad& e^b_{i,j}&\mapsto & x^b_{i,j},
\end{array}
$$
where $x^b_{i,j}=x_0^{j-i+\sum b_i}x_1^{-b_1}\dots x_n^{-b_n}$ is
a generator of $hom_{D^bCoh(\mathbb P^n)}(\mathcal O(i),\mathcal
O(j))$.

The next step is to show that this identification of morphisms
respect composition maps. For $i<j<k$, morphism spaces
$hom_{D(\mathcal{\widetilde D}/\Gamma)}(\mathcal U(i),\mathcal
U(j))$ and $hom_{D(\mathcal{\widetilde D}/\Gamma)}(\mathcal
U(j),\mathcal U(k))$ are nontrivial. The composition map is
computed as wedge product of cohomology in the following
decomposition.
\begin{eqnarray*}
&hom_{D(\mathcal {\widetilde D}/\Gamma)}(\mathcal U(j),\mathcal
U(k))\otimes hom_{D(\mathcal {\widetilde D}/\Gamma)}(\mathcal
U(i),\mathcal U(j))\\
=&\bigoplus\limits_{\mathcal U^b(j)\subset \mathcal U^0(i)}\quad
\bigoplus\limits_{\mathcal U^{b+c}(k)\subset \mathcal U^b(j)}H^0
(\Omega(\overline{\mathcal U^{b+c}(k)},d))\otimes
H^0(\Omega(\overline{\mathcal U^b(j)},d)).
\end{eqnarray*}
Therefore one observes that the multiplication restricted on
subspaces of subspaces $ H^0(\Omega(\overline{\mathcal
U^{b+c}(k)},d))$ and $H^0(\Omega(\overline{\mathcal
U^{b+c}(k)},d))$ gives rise to $e^c_{j,k}\cdot
e^b_{i,j}=e^{b+c}_{i,k}$.

\begin{center}
\includegraphics{fig.5}

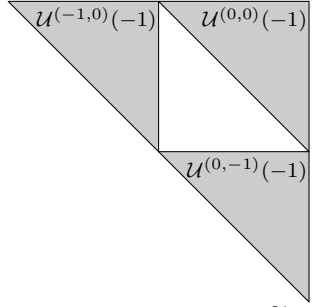
\parbox{8cm}{Fig.4 The case of $\mathbb P^2$. On $(\widetilde S^1)^2$, the base of
$\widetilde M$, $\mathcal U^{(0,0)}(-1)$, $\mathcal U^{(0,1)}(-1)$
and $\mathcal U^{(1,0)}(-1)$ (three shaded triangles) inside
$\mathcal U^{(0,0)}(-2)$ (the larger triangle) correspond to three
generators of $hom_{D^b Coh (\mathbb P^2)}(\mathcal O(-2),\mathcal
O(-1))$.}
\end{center}

Under the isomorphism $\nu$, $e^c_{j,k}$ and $e^b_{i,j}$
correspond to $x^c_{j,k}$ and $x^b_{i,j}$ respectively. The
composition
$$\begin{array}{lll}x^c_{j,k}\cdot x^b_{i,j}& = & x_0^{k-j+\sum c_i}x_1^{-c_1}\dots
x_n^{-c_n}\cdot x_0^{j-i+\sum b_i}x_1^{-b_1}\dots x_n^{-b_n}\\ \
&=&x_0^{k-i+\sum (b_i+c_i)}x_1^{-b_1-c_1}\dots x_n^{-b_n-c_n}\\
\ &=&x^{b+c}_{i,k},\end{array}$$which means $\nu(e^c_{j,k}\cdot
e^b_{i,j})=x^c_{j,k}\cdot x^b_{i,j}$. Therefore we have shown the
quiver structures of two full strong exceptional collections are
the same, the equivalence of triangulated categories $D^bCoh(\mathbb
P^n)\cong D(\widetilde {\mathcal D}/\Gamma)$ follows. By Lemma
\ref{passtosheaves}, our main theorem holds.
\begin{theo*}
$D \mathcal{ F}$ is equivalent to $D^b Coh(\mathbb P^n)$.\qed
\end{theo*}

\begin{rema}
The Lagrangians $\mathcal L^a(i)$ for $i\in \{-n-1,\dots, -1\}$
are standard branes over particular open sets (after the lift with
respect to the covering map $\pi$). They correspond to certain
constructible sheaves with respect to a fixed stratification
$\widetilde \Lambda$. Pushing forward to the constructible sheaves
on $(S^1)^n$ via the covering map $\pi$, they are constructible
sheaves on $(S^1)^n$ with respect to the stratification $\Lambda$,
which coincides with the stratification given by Bondal \cite{b}.
Bondal arrives at this stratification from the coherent sheaves on
a toric variety, while we obtain this stratification via the
standard Lagrangians constructed by T-duality.
\end{rema}


\begin{thebibliography}{CK}
\bibitem{ab} M. Abouzaid, \emph{Homogeneous coordinate rings
and mirror symmetry for toric varieties}, Geom. Topol. 10 (2006),
1097-1157 (math.SG/0511644). \
\bibitem{ab2} M. Abouzaid, \emph{Morse homology, tropical
geometry, and homological mirror symmetry for toric varieties},
math.SG/0610004.
\bibitem{ap} D. Arinkin, A. Polishchuk, \emph{Fukaya category and Fourier transform}, Winter School on Mirror Symmetry, Vector Bundles and Lagrangian Submanifolds (Cambridge, MA, 1999), 261 - 274,
AMS/IP Stud. Adv. Math., 23, Amer. Math. Soc., Providence, RI,
2001 (math/9811023).
\bibitem{au} D. Auroux, \emph{Mirror symmetry and T-duality in the complement of an anticanonical
divisor}, arXiv: 0706.3207
\bibitem{ako} D. Auroux, L. Katzarkov, D. Orlov, \emph{Mirror
symmetry for weighted projective planes and their noncommutative
deformations}, math.AG/0404281
\bibitem{ako2} D. Auroux, L. Karzarkov, D. Orlov, \emph{Mirror
symmetry for Del Pezzo surfaces: vanishing cycles and coherent
sheaves}, Inventiones Math. 166 (2006), 537-582 (math.AG/0506166).
\bibitem{b} A. Bondal, \emph{Integrable systems related to
triangulated categories}, talk at MSRI workshop on Generalized
McKay Correspondecnes and Representation Theory, regarding the
joint work with W.-D. Ruan, 3/21/2006 (video is at
www.msri.org/communications/vmath/VMathVideos/VideoInfo/2469\\/
show\underline \ video).
\bibitem{b2} A. Bondal, \emph{Derived categories of toric varieties}, Convex and Algebraic geometry, Oberwolfach
conference reports, EMS Publishing House 3 (2006), 284-286.
\bibitem{lc} K. Chan and N. C. Leung, \emph{Mirror symmetry for toric fano manifolds via SYZ
transformations}, arXiv: 0801:2830.
\bibitem{co} C.-H. Cho, Y.-G. Oh, \emph{Floer cohomology and disc instantons of Lagrangian torus fibers in Fano
toric manifolds}, Asian J. Math. 10 (2006), 773-814
(math.SG/0308225).
\bibitem{fltz} B. Fang, C-C. Liu, D. Treumann, E. Zaslow, \emph{T-duality and equivariant homological mirror symmetry for toric varieties}, in preparation.
\bibitem{fooo} K. Fukaya, Y.-G. Oh, H. Ohta, K. Ono,
\emph{Lagrangian intersection Floer Theory: anomaly and
obstruction}, preprint, second expanded version, 2006.
\bibitem{fss} K. Fukaya, P. Seidel, I. Smith, \emph{Exact Lagrangian submanifolds in simply-connected cotangent
bundles}, math/0701783.
\bibitem{h} K. Hori, \emph{Mirror symmetry and quantumn geometry},
Proceedings of the International Congress of Mathematicians, Vol.
III (Beijing, 2002), 431 - 443, Higher Ed. Press, Beijing, 2002
(hep-th/0207068).
\bibitem{hiv} K. Hori, A. Iqbal, C. Vafa, \emph{D-branes and
mirror symmetry}, hep-th/0005247.
\bibitem{hv} K. Hori, C. Vafa, \emph{Mirror symmetry},
hep-th/0002222.
\bibitem{ke} B. Keller, \emph{On the cyclic homology of exact
categories}, J. Pure Appl. Algebra 136(1999), no. 1, 1-56.
\bibitem{k} M. Kontsevich, \emph{Lectures at ENS, Paris, Spring 1998}, notes taken by
J. Bellaiche, J.-F. Dat, I. Marin, G. Racinet and H.
Randriambololona.
\bibitem{k2} M. Kontsevich, \emph{Homological algebra of mirror
symmetry}, Proc. International Congress of Mathematicians
(Z\"urich, 1994), Birkh\"auser, Basel, 1995, 120-139.
\bibitem{lv}
N. C. Leung, C. Vafa, \emph{Branes and toric geometry}, Adv.
Theor. Math. Phys. 2 (1998), no. 1, 91 - 118.
\bibitem{lyz} N. C. Leung, S-T Yau, E. Zaslow, \emph{From Special Lagrangian to Hermitian-Yang-Mills
via Fourier-Mukai Transform}, Adv. Theor. Math. Phys. 4 (2000),
no. 6, 1319-1341.
\bibitem{n} D. Nadler, \emph{Microlocal branes are constructible
sheaves}, math/0612399.
\bibitem{nz} D. Nadler, E. Zaslow: \emph{Constructible
sheaves and the Fukaya Category}, math/0604379.
\bibitem{syz} A. Strominger, S.-T. Yau, E. Zaslow, \emph{Mirror symmetry
is T-duality}, Nucl. Phys. B 479 (1996), 243-259 (hep-th/9606040).
\end{thebibliography}
\end{document}